\documentclass[12pt,twoside]{article}

\usepackage{amsmath,amssymb}
\usepackage{amsthm}
\usepackage{stmaryrd}
\usepackage{fancyhdr}
\usepackage{times}
\usepackage{mathrsfs} % Required for C and H
\usepackage{sectsty}
\usepackage[pdftex,dvips]{geometry}

\input diagxy.tex

\sectionfont{\small\normalfont\bfseries\scshape\centering\MakeUppercase }

\subsectionfont{\small\normalfont\bfseries\scshape }

\makeatletter
  \def\@seccntformat#1{\csname the#1\endcsname.\quad}
\makeatother

\newtheoremstyle{fact}% name
     {\topsep}%      Space above
     {\topsep}%      Space below
     {\slshape}%         Body font
     {}%         Indent amount (empty = no indent, \parindent = para indent)
     {\bfseries}% Thm head font 
     {}%        Punctuation after thm head
     { }%     Space after thm head: " " = normal interword space;
     {\thmname{#1}\thmnumber{ #2.}\thmnote{ \rm (#3)}}
\newcommand{\cat}[1]{\mathfrak{#1}}

\newcommand{\rest}{\mbox{\parbox[t]{0.1cm}{$|$\\[-10pt] $|$}}}

\newcommand{\Invisible}[1]{}

\newbox\pbbox
\setbox\pbbox=\hbox{\xy \POS(75,0)\ar@{-} (0,0) \ar@{-} (75,75)\endxy}

\newtheorem{theorem}{Theorem}[section]
\newtheorem*{theorem*}{Theorem} 
\newtheorem{lemma}[theorem]{Lemma}
\newtheorem{proposition}[theorem]{Proposition}
\newtheorem{corollary}[theorem]{Corollary}

\theoremstyle{definition}
\newtheorem{remark}[theorem]{Remark}
\newtheorem*{remark*}{Remark}

\newtheorem*{question*}{Question}
\newtheorem*{examples*}{Examples}  
\newtheorem{example}[theorem]{Example}
\newtheorem*{example*}{Example}

\theoremstyle{fact}
\newtheorem{fact}{Fact}

\def\proofont{\fontseries{bx}\fontshape{sc}\selectfont}
\def\proofname{Proof. }

\newcommand{\pcite}[2]{{\cite[#1]{#2}}}

\makeatletter
\renewenvironment{proof}[1][\proofname]{\par
  \normalfont
  \topsep6\p@\@plus6\p@ \trivlist
  \item[\hskip\labelsep\noindent\proofont #1]\ignorespaces
}{%
  \qed\endtrivlist
}
\makeatother

\author{G\'abor Luk\'acs
\thanks{I gratefully acknowledge the generous financial support received
from the Alexander von Humboldt Foundation, the Killam Trusts, and 
Dalhousie University that enabled me to do this research.}}
   
\title{Lifted closure operators
\thanks{2000 Mathematics Subject Classification: 18A32 54B30 
(18B30 22A05 46K10)}}

\begin{document}

\makeatletter
\let\mytitle\@title
\chead{\small\itshape G. Luk\'acs / \mytitle }
\fancyhead[RO,LE]{\small \thepage}
\makeatother

\maketitle

\def\thanks#1{} 

\thispagestyle{empty}

%Important, but not explicitly quoted references
\nocite{Atlantis}

\begin{abstract}
In this paper, we study the properties of closure operators obtained as 
initial lifts along a reflector, and compactness with respect to them in 
particular.  Applications in the areas of topology, topological groups and 
topological $*$-algebras are presented.
\end{abstract}

\section{Introduction}

Throughout this paper, $\cat{X}$ and $\cat{A}$ are finitely-complete
categories with proper $(\mathcal{E},\mathcal{M})$ and
$(\mathcal{F},\mathcal{N})$ factorization systems, respectively
(cf.\ \cite[7.2]{HerrRC}). For $X \in \cat{X}$,
the class $\mathcal{M}/X$ of subobjects of $X$ 
is preordered by the relation $m \leq n$ $\Longleftrightarrow$ 
$(\exists j)\, m=nj$, and denoted by $\operatorname{sub} X$.

A {\em closure operator} with respect to $\mathcal{M}$ is a family
$c=(c_X)_{X \in \cat{X}}$ of maps 
$c_X\colon \operatorname{sub} X \rightarrow \operatorname{sub} X$ such 
that $m \leq c_X(m)$, $c_X(m) \leq c_X(n)$ if $m \leq n$, and 
$f(c_X(m)) \leq c_Y(f(m))$  for all $f\colon X \rightarrow Y$ in 
$\cat{X}$ and $m,n \in \operatorname{sub} X$ (cf.\ \cite{DT}). The term 
``closure operator" has been used for almost a century in various
meanings in the context of categories of topological spaces and lattices.
However, this (general) categorical notion of closure operators was 
invented by Dikranjan and Giuli \cite{DikGiuCO}, and was further developed 
by Dikranjan, Giuli and Tholen in \cite{DGT}. Following \cite{DikGiu}, 
an object $X \in \cat{X}$ is said to be {\em $c$-compact} if for every 
$Y\in\cat{X}$, the projection $\pi_Y \colon X \times Y \rightarrow Y$ is 
{\em $c$-preserving}, or in other words, if $\pi_Y(c_X(m))=c_Y(\pi_Y(m))$
for every $m \in \operatorname{sub} (X\times Y)$. 

If $F \dashv U \colon \cat{A} \longrightarrow \cat{X}$ is an adjunction 
with unit $\eta$, and $c$ is a closure operator of $\cat{A}$ with respect to 
$\mathcal{N}$, then one defines the {\em $F$-initial lift of $c$} as
\begin{equation} \label{eq:init}
c_X^\eta (m) := \eta_X^{-1}(U c_{FX}(Fm))
\end{equation}
for every $X \in \cat{X}, m \in \operatorname{sub}(X)$ (cf.\ \cite[5.13]{DT}). 
The morphism $Fm$ need not belong to $\mathcal{N}$, so in order to ensure 
that (\ref{eq:init}) is defined, one extends the  notion of closure as 
$c_{FX}(Fm) = c_{FX}((Fm)(1_{FM}))$ (cf.\ \cite[5.7]{DT}).

In this paper, we study the special case of this construction where
$\cat{A}$ is a reflective subcategory of $\cat{X}$ with reflector $F=R$
and reflection $\eta=\rho$. We provide a characterization of
$c^\rho$-compact objects in $\cat{X}$ under certain conditions
(Theorem~\ref{thm:compact:main}), and applications in the areas of 
topology, topological groups and topological $*$-algebras are presented
We also investigate the categorical properties
of $c$ that are inherited by $c^\rho$.

\pagebreak[2]

\bigskip

The paper is structured as follows: We start off with an investigation of
the properties of lifted closure operators in section~\ref{sect:basic},
with an emphasis on compactness with respect to them. 
In section~\ref{sect:tych}, we show that the notion of
w-compactness (introduced by T.~Ishii in \cite{Ishii3}) is equivalent to
being compact with respect to the initial lift of the Kuratowski closure
operator on $\mathsf{Tych}$ (Tychonoff spaces and their continuous maps)
along the Tychonoff functor 
$\tau\colon\mathsf{Top}\longrightarrow\nolinebreak\mathsf{Tych}$.
In section~\ref{sect:bohr}, the Bohr-closure $c_b$ for topological groups 
is presented as the initial lift of the Kuratowski closure along the
Bohr-compactification, and Theorem~\ref{thm:compact:main} is used to
characterize $c_b$-compactness. The Bohr-closure also turns out to be a 
convenient example to show that certain properties are not preserved by 
the initial lift along a reflection. Finally, in section~\ref{sect:reptop},
we prove that the $*$-representation topology on topological $*$-algebras
is an initial lift of the Kuratowski closure along a suitable 
reflection.

%The idea of using the concept of factorization systems in order to
%generalize topological properties, such as compactness, to arbitrary
%categories already appeared in the papers of Manes \cite{Manes},
%Herrlich, Salicrup, and Strecker \cite{Herr0}, and Castellini
%\cite{Cast}. Such generalizations, as well as ``generalized"  
%injectivity %($c$-injectivity), were successfully used for investigating
%modules by %Fay \cite{Fay0}, and Dikranjan and Giuli \cite{DikGiuMod},
%nilpotent, torsion, and locally nilpotent groups by Fay \cite{Fay4}, Fay
%and Walls \cite{Fay1}, \cite{Fay3}, \cite{Fay5}, and rings by Joubert

\section{Properties of lifted closure operators}
\label{sect:basic}

In this section, we present results of a positive nature, in other words,  
properties that carry over from $c$ to $c^\rho$. Counterexamples, 
which show that certain properties are not inherited by $c^\rho$, appear 
later on.

As stated in the Introduction, $\cat{A}$ is a reflective subcategory of 
$\cat{X}$ with reflector $R$ and reflection $\rho$, and thus 
(\ref{eq:init}) becomes
\begin{equation}
c_X^\rho(m) := \rho^{-1}_X(c_{RX}(Rm)).
\end{equation}
It can be shown that $\mathcal{N} \subseteq \mathcal{M}$ and  
$R \mathcal{E} \subseteq \mathcal{F}$ (cf.\ \cite[5.13]{DT}), so in order 
to avoid confusion, we put $\operatorname{sub}_\cat{A} A$ for the 
class of subobjects of $A \in \cat{A}$ in $\cat{A}$ (i.e., $\mathcal{N}/A$) 
when necessary.

It was observed by Baron \cite{Baron} that every reflection decomposes
into two epireflections.  Using the same idea, we factorize $R$ as an
$\mathcal{E}$-reflector followed by an epireflector: For
every $X \in \cat{X}$, one can factorize the reflection as
$\rho_X=\rho_X^\mathcal{M}\rho_X^\mathcal{E}$, where
$\rho_X^\mathcal{E}\in\mathcal{E}$ and $\rho_X^\mathcal{M}\in\mathcal{M}$,
and put $R^\mathcal{E} X$ for the codomain of $\rho_X^\mathcal{E}$.
Obviously, $R^\mathcal{E} X$ belongs to   
$\cat{B} = \{ S \in \cat{X} \mid \rho_S \in \mathcal{M} \}$, and
$\cat{B}$ is the $\mathcal{E}$-reflective hull  of $\cat{A}$ in 
$\cat{X}$ (naturally, with reflector $R^\mathcal{E}$ and reflection
$\rho^\mathcal{E}$). Furthermore, $\rho^\mathcal{M}$ is an epimorphism in
$\cat{B}$, so $\cat{A}$ is an epireflective subcategory of $\cat{B}$ with 
reflector $R^\mathcal{M} = R \rest_\cat{B}$, and
$R = R^\mathcal{M} R^\mathcal{E}$.
Since $\cat{B}$ is $\mathcal{E}$-reflective, it inherits from $\cat{X}$
a proper $(\mathcal{E}_\cat{B},\mathcal{M}_\cat{B})$-factorization system,
where $\mathcal{E}_\cat{B}= \mathcal{E} \cap \operatorname{mor} \cat{B}$ 
and $\mathcal{M}_\cat{B}=\mathcal{M}\cap\operatorname{mor}\cat{B}$. 

\begin{remark} \label{rem:crho-dense} 
One says that $m \in \operatorname{sub} X$ is {\em $c$-dense} if
$c_X(m)=1_X$. Every $S \in \cat{B}$ is $c^\rho$-dense in its reflection
$RS$, because $c_{RS}^\rho(\rho_S) = \rho^{-1}_{RS}(c_{RS}(R\rho_S)) =
c_{RS}(1_{RS}) = 1_{RS}$. Thus, for every $X \in \cat{X}$, $\rho_X$ has a
$c^\rho$-dense image in $RX$.
\end{remark}

\begin{proposition} \label{prop:inf}
Let $c$ be a closure operator on $\cat{A}$.
For every $X \in \cat{X}$ and $m \in \operatorname{sub} X$, the 
class $\mathcal{C}:= \{ \rho_X^{-1}(c_{RX}(n)) \mid 
m \leq \rho_X^{-1}(n), n \in \operatorname{sub}_\cat{A} \}$ has an 
infimum and $\bigwedge \mathcal{C} = c^\rho_X(m)$.
\end{proposition}

\begin{proof}
Clearly, $c^\rho_X(m) \in \mathcal{C}$, because $\rho_X(m)\leq (Rm)(1_{RM})$. 
On the other hand, if $m\leq\rho_X^{-1}(n)$ for 
$n\in\operatorname{sub}_\cat{A}$, then  $\rho_X(m) \leq n$, and thus 
$(Rm)(1_{RM}) \leq n$; in particular, $c_{RX}(Rm) \leq c_{RX}(n)$. 
Therefore, $c^\rho_X(m)$ is a lower bound to $\mathcal{C}$. 
\end{proof}

A closure operator $c$ is {\em idempotent} if
$c_A(c_A(m))=c_A(m)$ for every $A\in\cat{A}$ and 
$m\in\nolinebreak\operatorname{sub}_\cat{A} A$. 

\begin{proposition} \label{prop:idem}
Let $c$ be a closure operator of $\cat{A}$. If $c$ is idempotent, then so 
is $c^\rho$.
\end{proposition}

\begin{proof}
Let $X \in \cat{X}$ and $m \in \operatorname{sub} X$. 
Clearly, $(R\rho_X^{-1}(k))(1_{R(K \times_{RX} X)})\leq k$ for every 
\linebreak
$k\in\operatorname{sub}_\cat{A} RX$,
so for $k = c_{RX}(Rm)$, the image of $R \rho_X^{-1}(c_{RX}(Rm))$ is 
contained in $c_{RX}(Rm)$.
Since $c$ is idempotent, we obtain:
\begin{align}
c_X^{\rho}(c_X^{\rho}(m)) & = \rho_X^{-1}(c_{RX}(R \rho_X^{-1}(c_{RX}(Rm)))) \\
& \leq \rho_X^{-1}(c_{RX}(c_{RX}(Rm))) = \rho_X^{-1}(c_{RX}(Rm)))
= c_X^{\rho}(m).
\end{align}
Therefore, $c_X^{\rho}(m)=c_X^{\rho}(c_X^{\rho}(m))$, because $c_X^{\rho}$ 
is extensive.
\end{proof}

One says that a closure operator $c$ is 
{\em  {\rm (}finitely{\rm )} productive} if for every
(finite) family $\{A_i\}_{i\in I}$ of objects in $\cat{A}$ and 
their subobjects $m_i \in\nolinebreak\operatorname{sub}_\cat{A} A_i$, 
$c_{\prod\limits_{i\in I} A_i}(\prod m_i) = 
\prod\limits_{i\in I} c_{A_i}(m_i)$.

\begin{proposition} \label{prop:prod}
Suppose that $\mathcal{F}$ is {\rm (}finite{\rm )} productive, and 
let $c$ be a closure operator of $\cat{A}$.
If $c$ is {\rm (}finitely{\rm )} productive and $R$ preserves 
{\rm (}finite{\rm )} products, then $c^\rho$ is {\rm (}finitely{\rm )} 
productive.
\end{proposition}

\begin{proof}
Let $\{X_i\}_{i \in I}$ be a {\rm (}finite{\rm )}  family of objects in 
$\cat{X}$, and let $m_i \in \operatorname{sub} X_i$ be their subobjects. Put
$X = \prod\limits_{i\in I} X_i$ and $m = \prod\limits_{i \in I} m_i$.
Since $R$ preserves {\rm (}finite{\rm )} products, 
$\rho_X = \prod\limits_{i \in I} \rho_{X_i}$ and
$Rm = \prod\limits_{i \in I} Rm_i$, and one has 
$(Rm)(1_{RM}) = \prod\limits_{i \in I} (Rm_i) (1_{RM_i})$, because 
$\mathcal{F}$ is {\rm (}finitely{\rm )} productive. Therefore, we obtain
\begin{align}
c^\rho_{X}(m) & \stackrel{\text{def}}  =
\textstyle
\rho^{-1}_{X}(c_{RX}(\prod\limits_{i \in I} Rm_i)) = 
(\prod\limits_{i\in I}\rho_{X_i})^{-1}(\prod\limits_{i\in I}c_{RX_i}(Rm_i))\\
&= \textstyle \prod\limits_{i \in I} \rho_{X_i}^{-1} (c_{RX_i}(Rm_i)) 
 = \prod\limits_{i \in I} c_{X_i}^\rho(m_i),
\end{align}
because $c$ is {\rm (}finitely{\rm )} productive, and limits interchange 
with each other.
\end{proof}

For $A \in \cat{A}$, a subobject $m \in\operatorname{sub}_\cat{A} A$ is
{\em $c$-closed} if $c_A(m)=m$. One says that $A$ is {\em $c$-separated}
(or {\em $c$-Hausdorff}) if the diagonal $\delta_A \colon A \rightarrow A 
\times A$ is $c$-closed in $A$ (cf.\ \cite[4.1]{CGT}).

\begin{proposition} \label{prop:separ}
Let $c$ be a a closure operator of $\cat{A}$.
For $X \in \cat{X}$, suppose that {\rm (a)} $RX$ is $c$-separated and 
{\rm (b)} $\rho_{X\times X} = \rho_X \times \rho_X$. Then $X$ is 
$c^\rho$-separated if and only if $\rho_X$ is a monomorphism.
\end{proposition}

\begin{proof}
Using the assumptions, one can easily see that
\begin{equation}
c_{X\times X}^\rho(\delta_X) \stackrel{\text{def}} =
\rho^{-1}_{X\times X}(c_{R(X\times X)}(R\delta_X))\stackrel{\text{(a)}} = 
(\rho_X \times \rho_X)^{-1} (c_{RX\times RX}(\delta_{RX}))
\stackrel{\text{(b)}} = (\rho_X \times \rho_X)^{-1} (\delta_{RX}).
\end{equation}
This simple computation shows that $c_{X\times X}^\rho(\delta_X)$ is 
precisely the kernelpair of $\rho_X$.
\end{proof}

We turn to characterizing the $c^\rho$-compact objects, and to that end, 
until the end of this section, $\mathcal{E}$ is assumed to be pullback 
stable.

%\begin{lemma} \label{lemma:compact:basic}
%Let $d$ be a closure operator of $\cat{X}$. 
%If $X \in \cat{X}$ is $d$-compact, then:
%
%\begin{list}{{\rm (\alph{enumi})}}
%{\usecounter{enumi}\setlength{\labelwidth}{25pt}\setlength{\topsep}{2pt}
%\setlength{\itemsep}{0pt} \setlength{\leftmargin}{20pt}}
%
%\item
%every image of $X$ is $d$-compact;
%
%\item
%the image of $X$ is $d$-closed in every $d$-separated object.
%\end{list}
%\end{lemma}
%
%\begin{proof}
%(a) Since $\mathcal{E}$ is pullback stable, if $e \colon X \rightarrow Z$ is 
%in $\mathcal{E}$, then $e_1:=e\times 1_Y\colon X\times Y\rightarrow Z\times Y$
%is also in $\mathcal{E}$, and so $e_1 (e_1^{-1}(n)) = n$ for every 
%$n \in \operatorname{sub} (Z \times Y)$. Thus, 
%$\pi_Y(e_1^{-1} (n)) = \pi_Y(n)$, and
%$\pi_Y(d_{Z \times Y}(n)) = \pi_Y(e_1^{-1}(d_{Z \times Y}(n))) \geq
%\pi_Y(d_{X \times Y}(e_1^{-1}(n)))  = d_Y (\pi_Y(e_1^{-1}(n))) =
%d_Y(\pi_Y(n))$.
%The reverse relation holds for any morphism in $\cat{X}$, and 
%therefore the statement follows.
%
%(b) Let $f\colon X \rightarrow Y$ be a morphism into a $d$-separated 
%object $Y \in \cat{X}$. By (a), $f(1_X)$ is $d$-compact, and so we may 
%assume that $f=m$ is a subobject of $Y$ . Since $Y$ is 
%$d$-separated, $\delta_Y$ is $d$-closed in $Y \times Y$, and thus 
%$(m\times 1_Y)^{-1}(\delta_Y)= \langle 1_X,m \rangle 
%\colon X \rightarrow X \times Y$ is $d$-closed. Therefore,
%$m=\pi_Y(\langle 1_X,m \rangle)$ is $d$-closed, because $X$ is 
%$d$-compact.
%\end{proof}

\begin{proposition} \label{prop:compact:basic}
Let $c$ be a closure operator on $\cat{A}$. 
If $X \in \cat{X}$ is $c^\rho$-compact, then:
\begin{list}{{\rm (\alph{enumi})}}
{\usecounter{enumi}\setlength{\labelwidth}{25pt}\setlength{\topsep}{2pt}
\setlength{\itemsep}{0pt} \setlength{\leftmargin}{20pt}}

\item
$R^\mathcal{E}X$ is $c^\rho$-compact;

\item
if $RX$ is $c$-separated, then $\rho_X \in \mathcal{E}$ and $RX$ is 
$c$-compact.
\end{list}
\end{proposition}

\begin{proof}
Since $\mathcal{E}$ is pullback stable, the image of any $c^\rho$-compact 
object under a morphism is again $c^\rho$-compact (cf.
\cite[5.2(3)]{CGT}), and thus (a)  follows. In order 
to show (b), notice that $R^\mathcal{E}X$ is $c^\rho$-closed in $RX$ (cf.
\cite[5.2(2)]{CGT}), because $RX$ is $c$-separated. On the other hand, 
$R^\mathcal{E}X$ is $c^\rho$-dense in $RX$ (Remark~\ref{rem:crho-dense}), 
and therefore $R^\mathcal{E}X = RX$.
\end{proof}

For $X\in\cat{X}$, we say that {\em $R$ preserves products with $X$} if 
$\rho_{X\times Y} = \rho_X \times \rho_Y$ for every $Y \in \cat{X}$.

\begin{theorem} \label{thm:compact:main}
Suppose that $\mathcal{F} \subset \mathcal{E}$.
Let $c$ be a closure operator of $\cat{A}$, and let 
$X \in \cat{X}$. 
If $R$ preserves products with $X$,
$\rho_X \in \mathcal{E}$, and $RX$ is $c$-compact, then 
$X$ is $c^\rho$-compact.
\end{theorem}

\begin{remark}
The condition of $\mathcal{F} \subset \mathcal{E}$ implies that
$\mathcal{F} \subseteq R \mathcal{E}$, and therefore
$\mathcal{F} = R \mathcal{E}$ (cf.\ \cite[5.13]{DT}). \linebreak[2] Its 
role in the theorem is to ensure that the factorization of a 
morphism in $\cat{A}$ coincides with its factorization in $\cat{X}$.
\end{remark}

In order to prove Theorem~\ref{thm:compact:main}, a technical lemma is 
required. It is an easy consequence of the {\em Beck-Chevalley Property}
(cf.~\cite[3.8]{Facets1}), and thus its straightforward proof is omitted:

\begin{lemma} \label{lemma:proj}
Let $X_1,X_2,Y_1,Y_2 \in \cat{X}$, and let $e\colon X_1 \rightarrow X_2$ 
and $f\colon Y_1 \rightarrow Y_2$ be morphisms. 
Put $g = e \times f\colon X_1 \times Y_1 \rightarrow X_2 \times Y_2$.
If $e \in \mathcal{E}$ {\rm (}and $\mathcal{E}$ is pullback stable{\rm )}, 
then for every $n \in \operatorname{sub}(X_2 \times Y_2)$, 
$\pi_{Y_1}(g^{-1}(n)) = f^{-1}(\pi_{Y_2}(n))$.
\qed
\end{lemma}

\begin{proof}[Proof of Theorem~\ref{thm:compact:main}.] 
Let $Y \in \cat{X}$. We  show that $\pi_Y\colon X\times Y \rightarrow Y$ 
is $c^\rho$-preserving. To that end, let $m\in\operatorname{sub}(X\times Y)$.
First, note that
\begin{equation} \label{eq:compute:im}
\pi_{RY}((Rm)(1_{RM})) = R(\pi_Y(m))(1_{RN}),
\end{equation}
where $\pi_Y(m)\colon N \rightarrow Y$. Indeed, $\pi_Y m = \pi_Y(m)e$, 
where $e \in \mathcal{E}$, so $R(\pi_Y m) = R(\pi_Y(m))Re$. Thus, 
$\pi_{RY}((Rm)(1_{RM})) = R(\pi_Y m)(1_{RM}) = 
R(\pi_Y(m))(Re(1_{RM})) = R(\pi_Y(m))(1_{RN})$, as
$Re\colon RM \rightarrow RN$ belongs to $R\mathcal{E} = \mathcal{F}$
and $R$ preserves products with $X$ (so $R\pi_Y = \pi_{RY}$).
By Lemma~\ref{lemma:proj}, 
$\pi_{Y}(\rho^{-1}_{X \times Y}(n))=\rho^{-1}_Y(\pi_{RY}(n))$ for every 
$n \in \operatorname{sub}(RX \times RY)$, because $\rho_X \in\mathcal{E}$ 
and $\rho_{X\times Y} = \rho_X \times \rho_Y$. Thus, 
putting $n =\nolinebreak c_{RX \times RY}(Rm)$ yields:
\begin{align*}
\pi_Y(c_{X\times Y}^\rho(m)) 
& \stackrel{\text{def}}{=}
\pi_{Y}(\rho^{-1}_{X \times Y}(c_{RX \times RY}(Rm))) = 
\rho^{-1}_Y(\pi_{RY}(c_{RX \times RY}(Rm))) 
\mbox{.} \\
\intertext{Since $RX$ is $c$-compact, one has
$c_{RY}(\pi_{RY}(k))= \pi_{RY}(c_{RX \times RY}(k))$
for every $k \in \operatorname{sub}(RX \times RY)$. By letting
$k = (Rm)(1_{RM})$ and then applying $\rho^{-1}_Y$ to both sides, one obtains}
\rho^{-1}_Y(\pi_{RY}(c_{RX \times RY}(Rm))) 
& =
\rho^{-1}_Y(c_{RY}(\pi_{RY}((Rm)(1_{RM}) )))  \\
& \stackrel{\text{(\ref{eq:compute:im})}} =
\rho^{-1}_Y(c_{RY}(R(\pi_Y(m))(1_{RN})))
\stackrel{\text{def}}{=}
c^\rho_Y(\pi_Y(m)).
\end{align*}
Therefore, $c^\rho_Y(\pi_Y(m)) =\pi_Y(c^\rho_{X\times Y}(m))$, as desired.
\end{proof}

\begin{corollary} \label{cor:compact:Eref}
Suppose that $\mathcal{N} = \mathcal{M} \cap \operatorname{mor}\cat{A}$ 
and that $\cat{A}$ is $\mathcal{E}$-reflective in $\cat{X}$. Let $c$ be a 
closure operator of $\cat{A}$, and let $X \in \cat{X}$. If $R$ preserves 
products with $X$, then $X$ is $c^\rho$-compact if and only if 
$RX$ is $c$-compact.
\end{corollary}

\begin{proof}
Since $\cat{A}$ is $\mathcal{E}$-reflective in $\cat{X}$, the 
$(\mathcal{E},\mathcal{M})$-factorization system in $\cat{X}$ restricts to 
a proper $(\mathcal{E}_\cat{A}, \mathcal{M}_\cat{A})$-factorization system 
in $\cat{A}$, where 
$\mathcal{E}_\cat{A}=\mathcal{E}\cap\operatorname{mor}\cat{A}$ and 
$\mathcal{M}_\cat{A}=\mathcal{M}\cap\operatorname{mor}\cat{A}$. Because of 
the ``orthogonality" relation, $\mathcal{N}$ and $\mathcal{F}$ determine 
each other, and thus $\mathcal{N} = \mathcal{M}_\cat{A}$ implies
$\mathcal{F} = \mathcal{E}_\cat{A}$. So, the conditions of 
Theorem~\ref{thm:compact:main} are fulfilled, and therefore if $RX$ is 
$c$-compact, then $X$ is $c^\rho$-compact. On the other hand, 
by Proposition~\ref{prop:compact:basic}(a), if $X$ is $c^\rho$-compact, 
then so is $R^\mathcal{E} X = RX$; hence, $RX$ is $c$-compact, as desired.
\end{proof}

A combination of Theorem~\ref{thm:compact:main} and 
Proposition~\ref{prop:compact:basic} yields:

\begin{corollary} \label{cor:compact:csep}
Let $c$ be a closure operator of $\cat{A}$. Suppose that 
$\mathcal{F} \subset \mathcal{E}$ and every object in  $\cat{A}$ is 
$c$-separated. If $X \in \mathcal{X}$ is so that $R$ preserves products 
with $X$, then  $X$ is $c^\rho$-compact if and only if $\rho_X\in\mathcal{E}$ 
and $RX$ is $c$-compact.
\qed
\end{corollary}

The {\em discrete} closure operator $s$ is defined as $s_A(m)=m$ for every 
$m \in \operatorname{sub}_\cat{A} A$ and $A \in \cat{A}$.

\begin{corollary} \label{cor:compact:trivial}
Suppose that $\mathcal{F} \subset \mathcal{E}$, and let $X \in \cat{X}$. 
If $R$ preserves products with $X$, then  $X$ is $s^\rho$-compact if 
and only if $\rho_X \in \mathcal{E}$.
\qed
\end{corollary}

\section{Application I: The Tychonoff functor}
\label{sect:tych}

Put $\mathbb{I}$ for the closed unit interval. For any $X \in \mathsf{Top}$
(the category of topological spaces and their continuous maps), the 
evaluation map $\Phi_X\colon X \rightarrow \mathbb{I}^{C(X,\mathbb{I})}$,
mapping $x \in X$ to $(f(x))_{f \in C(X,\mathbb{I})}$, is continuous.
One sets $\tau X = \Phi_X (X)$, and $\tau$ is called the {\em Tychonoff 
functor}\,: It is the reflector of $\mathsf{Tych}$ (the full subcategory 
formed by the Tychonoff spaces) in $\mathsf{Top}$.  

Put $\cat{X}=\mathsf{Top}$, $\cat{A}=\mathsf{Tych}$, 
$R=\tau$, and $\rho=\theta$.  The category $\mathsf{Top}$ admits a proper
$(\mathsf{Onto},\mathsf{Embed})$-factorization system, and the 
reflection $\theta_X\colon X\rightarrow\tau X$ is surjective (by definition) for
every $X\in\nolinebreak\mathsf{Top}$. Thus, by putting  
$\mathcal{N}=\mathcal{M}_\cat{A}$, we arrive at the setting described in 
Corollary~\ref{cor:compact:Eref}, because $\mathcal{E}=\mathsf{Onto}$ is 
obviously pullback stable. Therefore, for any closure operator $c$ 
of $\mathsf{Tych}$, if $\tau$ preserves products with 
$X\in\mathsf{Top}$, then $X$ is $c^\theta$-compact if and only if 
$\tau X$ is $c$-compact. Here, we are interested in the case where 
$c=K$, the Kuratowski-closure. Since the $K$-compact spaces in
$\mathsf{Tych}$ are precisely the compact ones (by the classical 
Kuratowski-Mr\'{o}wka Theorem, cf. \cite[3.1.16]{Engel6}), the 
problem of characterizing $K^\theta$-compact spaces in $\mathsf{Top}$ 
boils down to a question concerning preservation of products by $\tau$. 
Hence, we use the terminology and the results of T.~Ishii, who studied
very carefully and extensively the Tychonoff functor in \cite{Ishii2} and
\cite{Ishii3}, and in \cite{Ishii} gave a survey of the results
concerning it.

\bigskip

Let $X \in \mathsf{Top}$; $F \subset X$ is a {\em zero-set} of $X$ if 
there exists a continuous map $f\colon X \rightarrow \mathbb{I}$ such that 
$F = \{ x \in X \mid f(x) = 0\}$, and $G \subset X$ is a {\em cozero-set}
of $X$ if there exists a continuous map $g\colon X \rightarrow \mathbb{I}$ 
such that $G = \{ x \in X \mid g(x) > 0\}$. A subset $A$ of $X$ is 
{\em $\tau$-open} if $A$ is a union of cozero-sets of $X$. (The complement 
of a $\tau$-open set is called a $\tau$-closed set.) Following T.~Ishii, 
a space $X$ is {\em w-compact} if for any family  $\{P_\lambda\}$ of 
$\tau$-open subsets of $X$ with the finite intersection property, 
$\bigcap \bar P_\lambda \neq \emptyset$ holds; or equivalently, if 
for any collection $\{F_\alpha\}$ of closed sets of $X$ that is closed 
under finite intersections and such that each $F_\alpha$ contains a 
non-empty cozero-set of $X$, one has $\bigcap F_\alpha \neq \emptyset$ 
(cf. \cite[1.4]{Ishii3}). Although $\tau X$ is compact for every w-compact 
space $X$ (cf. \cite[2.1]{Ishii3}), there exists a non-w-compact, regular 
Hausdorff space $X$ such that $\tau X$ is compact (cf. 
\cite[2.3]{Ishii3}).

\begin{theorem} \label{thm:wc}
A space $X \in \mathsf{Top}$ is $K^\theta$-compact if and only if it is 
w-compact.
\end{theorem}

In order to prove Theorem~\ref{thm:wc}, we need two results of T.~Ishii.
Recall that a continuous map $f\colon X \rightarrow Y$ of topological 
spaces is a {\em $Z$-map} if $f(Z)$ is closed in Y for every zero-set $Z$ 
of $X$.

\begin{samepage}
\begin{fact}[\pcite{2.7}{Ishii3}, \pcite{3.11}{Ishii}] \label{fact:wcZ}
For a space $X \in \mathsf{Top}$, the following properties are equivalent:

\begin{list}{{\rm (\roman{enumi})}}
{\usecounter{enumi}\setlength{\labelwidth}{30pt}\setlength{\topsep}{-0pt}
\setlength{\itemsep}{-2.5pt}\setlength{\labelsep}{6pt}}

\item
$X$ is w-compact;

\item
The projection $\pi_Y\colon X \times Y \rightarrow Y$ is a $Z$-map for 
every $Y \in \mathsf{Top}$.
\end{list}
\end{fact}
\end{samepage}

\begin{samepage}
\begin{fact}[\pcite{1.5}{Ishii3}, \pcite{4.1}{Ishii}] \label{fact:LwcPro}
For a space $X \in \mathsf{Top}$, the following statements are equivalent:

\begin{list}{{\rm (\roman{enumi})}}
{\usecounter{enumi}\setlength{\labelwidth}{30pt}\setlength{\topsep}{0pt}
\setlength{\itemsep}{-2.5pt}\setlength{\labelsep}{6pt}}

\item
for every $x \in X$ there exists a cozero-set $W$ containing $x$ such that
$\bar W$ is w-compact;

\item
$\tau(X \times Y) = \tau (X) \times \tau (Y)$ for every $Y \in 
\mathsf{Top}$.

\end{list}
\end{fact}
\end{samepage}

\begin{proof}[Proof of Theorem~\ref{thm:wc}.]
Suppose that $X$ is $K^\theta$-compact, and let $Y \in \mathsf{Top}$. 
In order to show that $\pi_Y\colon X \times Y \rightarrow Y$ is a $Z$-map, 
let $Z \subset X \times Y$ be a zero-set in $X \times Y$. In particular, 
$Z$ is $K^\theta$-closed in $X \times Y$, and thus $\pi_Y(Z)$ is 
$K^\theta$-closed in $Y$, because $X$ is $K^\theta$-compact and 
$K^\theta$ idempotent (Proposition~\ref{prop:idem}).
Therefore, $\pi_Y(Z)$ is closed in $Y$, which shows that $\pi_Y$ is a 
$Z$-map. This holds for every $Y\in\mathsf{Top}$, 
and hence $X$ is w-compact (Fact~\ref{fact:wcZ}).

Conversely, suppose that $X$ is w-compact. Then $\tau$ preserves products 
with $X$ (Fact~\ref{fact:LwcPro}, with $W=X$ in (i)). Obviously, 
$\mathcal{E}=\mathsf{Onto}$ is pullback stable, and thus, by
Corollary~\ref{cor:compact:Eref}, $X$ is $K^\theta$-compact if and only if
$\tau X$ is $K$-compact. The latter is equivalent to $\tau X$ being 
compact, by the classical Kuratowski-Mr\'{o}wka Theorem 
(cf. \cite[3.1.16]{Engel6}). Since $X$ is $w$-compact, $\tau X$ is compact 
(cf. \cite[2.1]{Ishii3}), and therefore $X$ is $K^\theta$-compact.
\end{proof}

We note that if $\tau$ preserves products with $X$ (which can be
guaranteed by the condition in Fact~\ref{fact:LwcPro}), then $\tau X$
being compact implies that $X$ is $K^\theta$-compact. Therefore, 
Theorem~\ref{thm:wc} has the following corollary.

\begin{corollary} Let $X \in \mathsf{Top}$ be such that for every $x \in 
X$ there exists a cozero-set $W$ containing $x$ such that $\bar W$ is 
w-compact. Then the following statements are equivalent:

\begin{list}{{\rm (\roman{enumi})}}
{\usecounter{enumi}\setlength{\labelwidth}{30pt}\setlength{\topsep}{-0pt}
\setlength{\itemsep}{-2.5pt}\setlength{\labelsep}{6pt}}

\item
$X$ is $K^\theta$-compact;

\item
$X$ is w-compact;

\item
$\tau X$ is compact.
\qed
\end{list}
\end{corollary}

\section{Application II: The Bohr-closure} 
\label{sect:bohr}

Similarly to the relationship between $\mathsf{Top}$ and 
$\mathsf{HComp}$ (compact Hausdorff spaces), the full subcategory
$\mathsf{Grp(HComp)}$ of compact Hausdorff groups is reflective in
$\mathsf{Grp(Top)}$ (topological groups and their continuous homomorphisms).
The reflection is called the {\em Bohr-compactification}, and is denoted 
by $\rho_G\colon G \rightarrow bG$; $\rho_G$ has a dense image, but it need 
not be onto. The kernel of $\rho_G$ is denoted by $\mathbf{n}(G)$ and it is 
called the {\em von Neumann kernel} of $G$.
There are many groups that do not admit any non-trivial 
continuous homomorphism into a compact group at all. Such groups are 
called {\em minimally almost periodic}, and their Bohr-compactification is 
trivial (cf. \cite{NeuWig}). Therefore, one should not expect $\rho_G$ to 
be an injection. Those groups $G$ for which $\mathbf{n}(G)=\{e\}$ are 
called {\em maximally almost periodic}. 

A group $G$ is {\em precompact} if for every neighborhood $U$ of $e\in G$,
there exists a finite subset $F \subset G$ such that $G=FU$.  It is 
interesting to note that while topological spaces $X$ whose Stone-\v{C}ech 
reflection $X\rightarrow\beta X$ is an embedding are exactly the Tychonoff 
spaces, groups $G$ such that $\rho_G$ is an embedding are characterized by 
being Hausdorff (and thus Tychonoff) and precompact.

\begin{remark} \label{rem:nG:ab}
It follows from the famous Peter-Weyl Theorem 
that finite-dimensional irreducible unitary 
representations of a compact Hausdorff group separate its points 
(cf.~\cite[Thm.~33]{Pontr}). Every such representation of an abelian 
group is one dimensional. Thus, for every compact Hausdorff abelian 
group $K$ and $k\in K$ such that $k\neq 0$, there is $\chi \in \hat K$ 
such that $\chi(k)\neq 0$, where \linebreak
$\hat K = \mathscr{H}(K,\mathbb{T})$ is the 
group of {\em continuous characters} of $K$ (i.e., continuous homomorphisms
$\chi\colon K\rightarrow\mathbb{T}$,
where $\mathbb{T}=\mathbb{R}/\mathbb{Z}$). 
If $G$ is an abelian group, then so is its Bohr-compactification $bG$ (because 
$\rho_G(G)$ is an abelian dense subgroup of $bG$). Therefore, $g \in\ker\rho_G$ 
if and only if $\chi(\rho_G(g))=0$ for every $\chi\in\widehat{bG}$. By the 
universal property of $bG$,  there is a one-to-one correspondence between 
$\hat G$ and $\widehat{bG}$, because $\mathbb{T}\in\mathsf{Grp(HComp)}$ - 
in other words, $\hat G = \widehat{bG}$ as sets. 
Hence, if $G$ is an abelian group, then 
\begin{equation} \label{eq:nG:ab}
\mathbf{n}(G)=\bigcap\limits_{\chi\in\hat G} \ker \chi.
\end{equation}
\end{remark}

The {\em Bohr-closure} $c_b$ of a subgroup $S \subseteq G$ is defined as 
$c_b(S) = \rho_G^{-1} ( \overline{\rho_G(S)})$. We equip $\mathsf{Grp(Top)}$
with the $(\mathsf{Onto},\mathsf{Embed})$-factorization system, while 
$\mathsf{Grp(HComp)}$ is provided with the 
$(\mathsf{Onto},\mathsf{ClEmb})$-factorization system ($\mathsf{ClEmb}$ 
stands for the closed embeddings). In this setting, it follows from 
Proposition~\ref{prop:inf} that $c_b$ is precisely the $b$-initial lift 
$s^\rho$ of the discrete closure operator $s$ on $\mathsf{Grp(HComp)}$.

\begin{proposition}
\mbox{ }

\begin{list}{{\rm (\alph{enumi})}}
{\usecounter{enumi}\setlength{\labelwidth}{25pt}\setlength{\topsep}{2pt}
\setlength{\itemsep}{0pt} \setlength{\leftmargin}{20pt}}

\item
$c_b$ is idempotent;

\item
$c_b$ is productive;

\item
$G \in \mathsf{Grp(Top)}$ is $c_b$-separated if and only if $\rho_G$ is 
injective {\rm (}i.e., $G$ is maximally almost periodic{\rm )};

\item
$G \in \mathsf{Grp(Top)}$ is $c_b$-compact if and only if $\rho_G$ is 
surjective.

\end{list}
\end{proposition}

\begin{proof}
(a) follows from Proposition~\ref{prop:idem}. (b) follows from 
Proposition~\ref{prop:prod}, because the Bohr-compactification preserves 
products, and surjective maps in $\mathsf{Grp(HComp)}$ are productive. (c) 
is an immediate consequence of Proposition~\ref{prop:separ}, because $b$ 
is productive and the class of monomorphisms coincides with the injective 
maps in $\mathsf{Grp(Top)}$. Finally, (d) is a special case of 
Corollary~\ref{cor:compact:trivial}.
\end{proof}

Let $G$ be a locally compact Hausdorff group. One says that $G$ is 
{\em central} if $G/Z(G)$ is compact (cf. \cite{GroMos1}, \cite{GroMos2}); 
$G$ is a {\em Moore group} if every irreducible unitary representation of 
$G$ is finite-dimensional (cf. \cite{Moore}, \cite{Robertson}). Moore 
groups are  automatically maximally almost periodic, because by the 
Gelfand-Ra\u{\i}kov Theorem (\cite[22.12]{HewRos}), irreducible 
representations of a locally compact group separate its points. 

\begin{theorem}
Let $G$ be a locally compact Hausdorff $c_b$-compact group.

\begin{list}{{\rm (\alph{enumi})}}
{\usecounter{enumi}\setlength{\labelwidth}{25pt}\setlength{\topsep}{2pt}
\setlength{\itemsep}{0pt} \setlength{\leftmargin}{20pt}}

\item
If $G$ is Moore, then it is compact.

\item
If $G$ is central, then it is compact.

\item
If $G$ has nilpotency class $\leq 2$, then it is compact.

\item
The quotient $G/\overline{[[G,G],G]}$ is compact; in particular,
$\mathbf{n}(G) \subseteq \overline{[[G,G],G]}$.

\item
Every continuous homomorphism $\alpha\colon G \rightarrow 
\langle \mathbb{R}, + \rangle$ is trivial, and thus 
$G$ is unimodular.
\end{list}
\end{theorem}

\begin{proof}
(a) Hughes \pcite{Theorem~1}{HughesPhD} \& \pcite{Theorem~1}{Hughes} 
showed that for a locally compact Hausdorff group $G$, a subspace 
$K \subset G$ is compact if and only if $K$ is compact in $G^w$, the group
$G$ equipped with the pointwise topology induced by its irreducible
unitary representations.  Remus and Trigos-Arrieta 
\cite[Theorem~1]{RemTri} showed that for Moore groups, $G^w = \rho_G(G)$. 
Since $G$ is $c_b$-compact, $\rho_G(G) = bG$, and therefore $G^w$ is 
compact. Hence, by Hughes' result, $G$ is compact.
(b) follows from (a), because every central group is Moore (cf. 
\cite[Theorem.~2.1]{GroMos2}).

(c) The locally compact abelian group $G_{ab}:=G/\overline{[G,G]}$, being the
continuous image of $G$, is $c_b$-compact. It is also Moore, because every
irreducible unitary representation of an abelian group is one dimensional.
Thus, by (a), $G_{ab}$ is compact. Since $G$ has nilpotency class $2$, one has 
$[[G,G],G]=1$, and so $\overline{[G,G]} \subseteq Z(G)$.
Therefore, $G/Z(G)$ is a quotient of $G_{ab}=G/\overline{[G,G]}$, and 
hence $G$ is central. The statement now follows by (b).

(d) follows from (c).

(e) $\alpha$ factors through $G_{ab}$, which is compact (by (d)), but the 
only compact subgroup of $\mathbb{R}$ is the trivial one. Therefore, 
$\alpha$  must be trivial.
\end{proof}

A group is {\em balanced} (or equivalently, {\em admits an invariant
basis}) if its left and right uniformities coincide. A Hausdorff locally
compact connected group is balanced if and only if it is central (cf.
\cite[Theorem~4.3]{GroMos1}), which yields:

\begin{corollary} \label{cor:cb}
Every Hausdorff locally compact connected balanced $c_b$-compact group is 
compact.
\end{corollary}

One might hope that $c_b$-compactness ``almost" implies compactness, but
this is not the case. For instance, every minimally almost periodic group 
has trivial Bohr-compactification, and in particular, they are $c_b$-compact.
The examples below show that all the conditions in Corollary~\ref{cor:cb}
are indeed necessary. 

\begin{example} \label{ex:SL2}
The special linear group $SL_2(\mathbb{R})$ is a locally compact, 
connected, and $c_b$-compact group, which is not compact.
Indeed, $SL_2(\mathbb{R})$ has the property that its
continuous surjective homomorphisms are open - in other words, it is {\em 
totally minimal} (cf. \cite[7.4]{DikProSto}, \cite{RemStoy}). In 
particular, every image of $SL_2(\mathbb{R})$ is complete, and thus
$c_b$-compact; in particular, 
$\rho_{SL_2(\mathbb{R})}(SL_2(\mathbb{R}))$ is a complete dense
subgroup of $b(SL_2(\mathbb{R}))$, and so $\rho_{SL_2(\mathbb{R})}$
is onto. Therefore, $b(SL_2(\mathbb{R}))$ is a quotient of 
$SL_2(\mathbb{R})$. The only non-trivial quotient of 
$SL_2(\mathbb{R})$ is $PSL_2(\mathbb{R})$, which is not compact, because 
the matrix $A =\left(\begin{array}{cc}1 & 1 \\ 0 & 1\end{array}\right)$ 
generates a discrete infinite subgroup. Hence, $SL_2(\mathbb{R})$ 
is minimally almost periodic.
\end{example}

\begin{example}
Nienhuys \cite{Nienhuys} showed that there exists a coarser metrizable 
group topology $\mathcal{T}$ on $\mathbb{R}$ such that 
$G:=(\mathbb{R},\mathcal{T})$ is minimally almost periodic. Thus,
$G$ is $c_b$-compact, and it is also balanced, because $G$ is abelian. 
Therefore, $G$ is a balanced, connected, and $c_b$-compact 
group, which is not compact.
\end{example}

Van der Waerden proved that every (algebraic)  homomorphism from a compact
connected semisimple Lie group into a compact group is continuous
(cf. \cite{vdW}). Groups $K$ with this property satisfy $K=\nolinebreak
b(K_d)$, where $K_d$ is the group $K$ equipped with the discrete topology,
and they are called {\em van der Waerden groups}, or briefly, {\em
vdW-groups} (cf.~\cite{HarKun}).  
For $D=K_d$, one obtains a discrete group that is $c_b$-compact.

\begin{example}
$SO_3(\mathbb{R})$ is a vdW-group, so $SO_3(\mathbb{R})_d$ is discrete,
infinite, and  $c_b$-compact. Thus, $SO_3(\mathbb{R})_d$ is a balanced, 
locally compact, and $c_b$-compact group, which is not compact. 
\end{example}

A closure operator $c$ is {\em hereditary} if for every 
$m\colon M \rightarrow A$ in $\operatorname{sub}_\cat{A} A$ and
$n\colon N \rightarrow N$ in $\operatorname{sub}_\cat{A} N$,
$c_M(n)=m^{-1}(c_A(n))$; $c$ is {\em weakly-hereditary} if for every 
$m\colon M \rightarrow A$ in $\operatorname{sub}_\cat{A} A$,
$m$ is $c$-dense in the domain of $c_A(m)$.

\begin{theorem}
The lift $c^\rho$ need not be weakly-hereditary when $c$ is hereditary, 
even if $\cat{A}$ is $\mathcal{E}$-reflective in $\cat{X}$  and 
$\mathcal{F}=\mathcal{E} \cap \operatorname{mor} \cat{A}$.
\end{theorem}

\begin{proof}
As we mentioned at the beginning of this section, groups whose 
Bohr-compactification is an embedding are precisely the precompact ones, 
and for them $c_b$ coincides with the Kuratowski closure. Hausdorff precompact 
groups form an $\mathsf{Onto}$-reflective subcategory of $\mathsf{Grp(Top)}$, 
and the reflection $G \rightarrow b^+G$ is given by the image of $G$ in $bG$.
Therefore, $c_b$ can also be seen as the $b^+$-initial lift of the 
hereditary Kuratowski closure operator of the 
precompact Hausdorff groups. Hence, it suffices to show that $c_b$ is not 
weakly-hereditary.

The Bohr-closure of the trivial subgroup of a topological group $G$
is simply its von Neumann kernel, $\mathbf{n}(G)$. Had $c_b$ been
weakly-hereditary, then 
\begin{equation}
\mathbf{n}(\mathbf{n}(G))=(c_b)_{\mathbf{n}(G)}(\{e\})=(c_b)_G(\{e\})
=\mathbf{n}(G)
\end{equation}
would have been true for every topological group $G$, but this is not
the case, as the Example~\ref{ex:from:Dik} shows.
\end{proof}

The example below was suggested by Dikran Dikranjan.

\begin{example} \label{ex:from:Dik}
For a prime number $p\neq 2$, let $A=Z(p^\infty)$, the Pr\"ufer group. 
It can be seen as the subgroup of $\mathbb{Q}/\mathbb{Z}$ generated
by the elements of $p$-power order, or the group formed by all  
$p^n$th roots of unity.
Since every proper subgroup of $A$ is finite, if $A$ is equipped
with a Hausdorff group topology in which it is not minimally almost 
periodic (in other words, $\mathbf{n}(A)\neq A$), then $\mathbf{n}(A)$ 
is a finite subgroup with the discrete topology, and thus 
$\mathbf{n}(\mathbf{n}(A))$ is trivial. Therefore, if $\tau$ is a 
Hausdorff group topology on $A$ such that $(A,\tau)$ is neither 
minimally nor maximally almost periodic, then 
$\mathbf{n}(\mathbf{n}(A,\tau))\neq \mathbf{n}(A,\tau)$.

Following Zelenyuk and Protasov \cite{ZelProt}, a sequence 
$\{a_n\}\subseteq G$ in a (discrete) abelian group $G$ is said to be 
a {\em $T$-sequence} if there is a Hausdorff group topology $\tau$ on $A$
such that $a_n\stackrel \tau \longrightarrow 0$. It follows from the proof of 
\cite[Theorem~1]{ZelProt} that if $\{a_n\}$ is a $T$-sequence
and $\tau$ is the finest group topology such that 
$a_n\stackrel \tau \longrightarrow 0$, then $(G,\tau)$ is universal in
the following sense: a homomorphism $\varphi\colon G \{a_n\} \rightarrow H$ 
into a topological group $H$ is continuous if and only if
$\varphi(a_n) \rightarrow 0$.

Let $c_n = \frac 1 {p^n}$ in $A$, put
\begin{equation}
b_n= -c_1 +c_{n^3 -n^2} + \cdots + c_{n^3-2n} + c_{n^3-n} +c_{n^3}
=-\frac 1 p + \frac 1 {p^{n^3 -n^2}}
+ \cdots + \frac 1 {p^{n^3-2n}} + \frac 1 {p^{n^3-n}} +\frac 1 {p^{n^3}},
\end{equation}
and consider the sequence $a_n$ defined as $b_1,c_1,b_2,c_2,b_3,c_3,\ldots$.
One can show that $\{a_n\}$ is a $T$-sequence in $A$; we chose to
omit the proof of this statement because of its technical nature and 
length. We set $\tau$ to be the finest Hausdorff group topology on $A$ 
such that $a_n \rightarrow 0$ in $\tau$, and apply (\ref{eq:nG:ab})
from Remark~\ref{rem:nG:ab} to show that 
$\mathbf{n}(A,\tau)=\langle \frac 1 p \rangle$.

To that end, let $\chi\in\widehat{(A,\tau)}$; then $\chi(a_n)\rightarrow 0$,
and in particular, $\chi(c_n) \rightarrow 0$. By \cite[Example~6]{ZelProt} 
and \cite[3.3]{DikMilTon}, $\chi(c_n) \rightarrow 0$ is equivalent to
$\chi=m\chi_1$, where $m\in\mathbb{Z}$ and $\chi_1$ is the natural embedding 
of $\mathbb{Z}(p^\infty)$ into 
$\mathbb{Q}/\mathbb{Z} \subseteq \mathbb{R}/\mathbb{Z}$. It is easily seen 
that $\chi_1(b_n) = - \frac 1 p$. Thus, if $\chi(b_n) \rightarrow 0$, then 
$-m \frac 1 p = 0$ in $\mathbb{R}/\mathbb{Z}$, and so $p \mid m$. On the other 
hand, it follows from $p\chi_1(b_n) \rightarrow -p\frac 1 p=0$ that 
$p\chi_1(a_n)\rightarrow 0$, so $p\chi_1$ is continuous with respect to 
$\tau$. Therefore,  $\widehat{(A,\tau)} = \{ m\chi_1 \mid m \in p \mathbb{Z}\}$. 
Hence, $\langle \frac 1 p\rangle \subseteq \ker \chi$ for
every $\chi \in \widehat{(A,\tau)}$. Since
$\ker p\chi_1 = \langle \frac 1 p\rangle$, this shows that
$\mathbf{n}(A,\tau) = \langle \frac 1 p \rangle$, as desired.
\end{example}

\section{Application III: The $\mathbf{*}$-representation topology}

\label{sect:reptop}

A {\em $*$-algebra} $A$ is an algebra over $\mathbb{C}$ with an
involution $\empty^*\colon A \rightarrow A$ such that
$(a+\lambda b)^* = a^*+\bar \lambda b^*$ and $(ab)^*=b^*a^*$
for every $a,b\in A$ and $\lambda \in \mathbb{C}$.
A {\em topological $*$-algebra} is a $*$-algebra $A$
and a topology on $A$ making the operations (addition, multiplication, 
additive inverse, involution) jointly continuous. 
The category of topological $*$-algebras
and their continuous $*$-homomorphisms is denoted by $\mathsf{T^*A}$. 
A {\em $C^*$-{\rm [}semi{\rm ]}norm} on a $*$-algebra $A$ is a [semi]norm 
$p$ that is submultiplicative and satisfies the $C^*$-identity - in other 
words, $p(ab)\leq p(a) p(b)$ and $p(a^* a)=p(a)^2$ for every $a,b\in A$.
We put $\mathcal{N}(A)$ for the set of continuous $C^*$-seminorms
on a topological $*$-algebra $A$. For every $p,q\in\mathcal{N}(A)$, 
one has $\max\{p,q\}\in\mathcal{N}(A)$ ($\max\{p,q\}$ is defined pointwise), 
which turns $\mathcal{N}(A)$ into a directed set.

A {\em $C^*$-algebra} is a complete Hausdorff topological algebra
whose topology is given by a single $C^*$-norm. The full subcategory of
$\mathsf{T^*A}$ formed by the $C^*$-algebras is denoted by $\mathsf{C^*A}$.
For $A \in \nolinebreak\mathsf{T^*A}$, a {\em $*$-representation} of $A$
is a continuous $*$-homomorphism $\pi\colon A\rightarrow B(H)$ 
of $A$ into the $C^*$-algebra of bounded operators on some Hilbert space $H$
(i.e., a morphism in $\mathsf{T^*A}$). The class
of $*$-representations of $A$ is denoted by $\mathcal{R}(A)$.

Let $p\in \mathcal{N}(A)$; 
$\ker p = \{ a\in A \mid p(a)=0\}$ is a $*$-ideal in $A$, and
$p$ induces a $C^*$-norm on the quotient $A/\ker p$, so 
the completion $A_p$ of this quotient with respect to 
$p$ is a $C^*$-algebra. By the celebrated Gelfand-Na\u{\i}mark-Segal theorem,
every $C^*$-algebra is $*$-isomorphic (and thus isometric) to a closed
subalgebra of $B(H)$ for a large enough Hilbert space $H$ 
(cf. \cite[2.6.1]{DixmierC}). Thus, we obtain a $*$-representation
$\pi_p\colon A \rightarrow A/\ker p \rightarrow A_p \rightarrow B(H)$ 
such that $p(x)=\|\pi_p(x)\|$. Conversely, each $\pi \in \mathcal{R}(A)$ 
gives rise to a $C^*$-seminorm $p_\pi(x)=\|\pi(x)\|$. Therefore, the initial
topology $\mathcal{T}_A$ induced by the class $\mathcal{R}(A)$ coincides 
with the one induces by the family of $C^*$-seminorms $\mathcal{N}(A)$.
The topology $\mathcal{T}_A$ is called the {\em $*$-representation topology},
and it defines a closure operator $c^*$ on subalgebras in $\mathsf{T^*A}$.
The closure $c^*_A(\{0\})$ of the trivial subalgebra is a closed $*$-ideal of 
$A$; it is called the {\em reducing ideal} of $A$, and 
denoted by $A_R$ (cf.~\cite[9.7]{Palmer2}). Notice that 
\begin{equation}
c_A^*(\{0\})=A_R=\bigcap\limits_{p \in \mathcal{N}(A)} \ker p = 
\bigcap\limits_{\pi \in \mathcal{R}(A)} \ker \pi,
\end{equation}
so $\mathcal{T}_A$ is Hausdorff (i.e.,  $A_R=0$) if and only if the 
$*$-representations of $A$ separate the points of $A$.

A topological $*$-algebra $A$ is said to be a {\em pro-$C^*$-algebra} 
if it is Hausdorff, complete, and its topology is generated by 
a family of $C^*$-seminorms - in other words, the topology of $A$
coincides with $\mathcal{T}_A$, it is complete, and $A_R=0$.
One can show that $A$ is a pro-$C^*$-algebra if and only if
$A$ is the limit in $\mathsf{T^*A}$  of $C^*$-algebras 
(cf.~\cite[1.1.1]{Phillips2}).
The full subcategory of $\mathsf{T^*A}$ consisting of pro-$C^*$-algebras 
is denoted by $\overline{\mathsf{P^*A}}$. These algebras were
studied under various names (locally $C^*$-algebras,
LMC*-algebras, etc.) by numerous authors; for instance, in section~3 of 
\cite{DubPor}, Dubuc and Porta essentially characterized commutative
pro-$C^*$-algebras (up to a $\mathsf k$-ification). Alluding only to a few 
more highlights, we mention the works of Inoue \cite{Inoue1},
Schm\"udgen \cite{Schmud}, Phillips \cite{Phillips1} \& \cite{Phillips2},
Bhatt and Karia \cite{BhatKar}, and Inoue and K\"ursten \cite{InoKur}.

\begin{proposition} \label{prop:TA-ref}
$\overline{\mathsf{P^*A}}$ is a reflective subcategory of \,$\mathsf{T^*A}$.
\end{proposition}

\begin{proof}
Let $A \in \mathsf{T^*A}$. For each $p\in\mathcal{N}(A)$, we put 
$A_p$ to be the completion of $A/\ker p$ with respect to the
$C^*$-norm that $p$ induces on it; obviously, $A_p$ is a $C^*$-algebra.
For every $q \geq p $ in the directed set $\mathcal{N}(A)$, there is a
surjective morphism $A_q \rightarrow A_p$. This defines a functor
\begin{align}
A_*\colon \mathcal{N}(A)&\longrightarrow\mathsf{C^*A}\subset\mathsf{T^*A}\\
p & \longmapsto A_p,
\intertext{and so we define the reflector} 
\overline{PC^*}(A) &= \lim A_*,
\end{align}
where the limit is taken in $\mathsf{T^*A}$. It follows from the definition 
that $\overline{PC^*}(A)$ is a pro-$C^*$-algebra.

To show that $\overline{PC^*}(-)\colon \mathsf{T^*A} \longrightarrow 
\overline{\mathsf{P^*A}}$ is a
functor, let $\varphi\colon A\rightarrow B$ be a morphism in $\mathsf{T^*A}$.
Then $\mathcal{N}(\varphi)\colon \mathcal{N}(B)\rightarrow \mathcal{N}(A)$
defined by $r \mapsto r  \varphi$ is an order-preserving map, and
$\varphi_r\colon A_{r \varphi} \rightarrow B_r$ is a natural transformation
from $A_* \mathcal{N}(\varphi)$ to $B_*$. Thus, there is a morphism
$\lim A_* \mathcal{N}(\varphi) \rightarrow \lim B_*$, and therefore
one sets
$\overline{PC^*}(\varphi)=(\lim A_*\rightarrow\lim A_*\mathcal{N}(\varphi)
\rightarrow\lim B_*)$.

To complete the proof, one defines 
$\eta_A\colon A \rightarrow \overline{PC^*}(A)$ by setting 
$\eta_A(x) = (x+\ker p)_{p \in \mathcal{N}(A)}$. Obviously,
$\eta_A$ is a morphism (because each $A \rightarrow A_p$ is so), and it is 
a natural transformation. 
Let $x=(x_p)_{p \in \mathcal{N}(A)}\in \overline{PC^*}(A)$, and pick 
$p_1,\ldots,p_k \in \mathcal{N}(A)$. Then 
$q=\max\{p_1,\ldots,p_k\}\in\mathcal{N}(A)$, and $A_q\rightarrow A_{p_i}$ is 
surjective. The image of $A$ under the canonical morphism $A\rightarrow A_q$
is dense, so for every $\varepsilon > 0$, 
there is $a \in A$ such that $q(\bar a_q-x_q)<\varepsilon$, and 
thus $p_i(\bar a_{p_i} - x_{p_i}) < \varepsilon$ (where $\bar a_p$ stands 
for the image  of $a$ in $A_p$). Therefore, the image $\eta_A(A)$ is dense in 
$\overline{PC^*}(A)$. 
If $B \in \overline{\mathsf{P^*A}}$, then 
$\eta_B=\operatorname{id}_B$. Hence, if $\varphi\colon A\rightarrow B$
is a morphism from $A\in\mathsf{T^*A}$ into $B\in\overline{\mathsf{P^*A}}$, then
$\varphi=\eta_B\varphi=\overline{PC^*}(\varphi)\eta_A$.
This decomposition is unique, because $\eta_A(A)$ is dense in 
$\overline{PC^*}(A)$.
\end{proof}

We turn now to fitting all these into the setting of section~\ref{sect:basic}.
Equip $\cat{X}=\mathsf{T^*A}$ with the usual
$(\mathsf{Onto},\mathsf{Embed})$-factorization system, and 
provide $\cat{A}=\overline{\mathsf{P^*A}}$ with the
$(\mathsf{Dense},\mathsf{ClEmb})$-factorization system
($\mathsf{Dense}$ and $\mathsf{ClEmb}$ stand for the maps
with a dense image and closed embeddings, respectively).
The $\mathsf{Onto}$-reflective hull of $\cat{A}$ in $\cat{X}$ is the 
category $\cat{B}=\mathsf{P^*A}$ consisting of the algebras $A\in\mathsf{T^*A}$
such that $A_R=0$ and $\mathcal{T}_A$ coincides with the topology of $A$;
the reflector $PC^*$ is given by $A \longmapsto \eta_A(A)$,
and the reflection $\alpha_A\colon A \rightarrow PC^*(A)$ is the
same as $\eta$, but with different codomain. The category $\cat{B}$ inherits 
the $(\mathsf{Onto},\mathsf{Embed})$-factorization system of $\cat{A}$.

\begin{theorem} \mbox{ } 

\begin{list}{{\rm (\alph{enumi})}}
{\usecounter{enumi}\setlength{\labelwidth}{25pt}\setlength{\topsep}{2pt}
\setlength{\itemsep}{0pt} \setlength{\leftmargin}{20pt}}

\item
$c^*=K^\alpha=K^\eta$, in other words, the $*$-representation 
topology is the $PC^*$-initial {\rm [}$\overline{PC^*}$-initial{\rm ]}
lift of the Kuratowski closure $K$ on $\mathsf{P^*A}$ 
{\rm [}$\overline{\mathsf{P^*A}}${\rm ]}.

\item
$c^*$ is finitely productive.

\item
An algebra $A \in \mathsf{T^*A}$ is $c^*$-separated if and only if
$\alpha_A$ is injective, or equivalently, $A_R=0$; in particular,
each algebra in $\mathsf{P^*A}$ is $c^*$-separated.

\item
An algebra $A \in \mathsf{T^*A}$ is $c^*$-compact in $\mathsf{T^*A}$ 
if and only if $PC^*(A)$ is $K$-compact in $\mathsf{P^*A}$; in this case,
$\eta_A$ is surjective.
\end{list}

\end{theorem}

\begin{proof}
(a) It follows from the definition that $\mathcal{T}_A$ is the 
initial topology on $A$ induced by the family of *-homomorphisms 
$(A \rightarrow\nolinebreak A_p)_{p\in\mathcal{N}(A)}$. Equivalently, 
$\mathcal{T}_A$ is the initial topology induced by the single *-homomorphism
$A \rightarrow\!\! \prod\limits_{p\in\mathcal{N}(A)}\!\!A_p$,
whose image is contained in $PC^*(A)\subseteq\overline{PC^*}(A)$.
Therefore, the statement follows from Proposition~\ref{prop:inf}.

(b) Since $PC^*(-)$ is additive, it preserves finite products in
$\mathsf{T^*A}$. Therefore, the statement follows from
Proposition~\ref{prop:prod}, as $K$ is productive,
and surjective morphisms in $\mathsf{P^*A}$ are productive. 

(c) follows from Proposition~\ref{prop:separ}, because $PC^*(-)$ preserves 
finite products, and monomorphism in $\mathsf{T^*A}$ are injections.

(d) It would be tempting to apply Theorem~\ref{thm:compact:main} to $\cat{A}$, 
but unfortunately, the condition of $\mathcal{F}\subset\nolinebreak\mathcal{E}$ 
fails: in our case, $\mathcal{F}=\mathsf{Dense}$ in $\overline{\mathsf{P^*A}}$ 
while $\mathcal{E}=\mathsf{Onto}$. Instead, we apply 
Corollary~\ref{cor:compact:Eref} to $\mathcal{B}=\nolinebreak\mathsf{P^*A}$ 
in order to obtain the first 
statement, while the second one is a consequence of 
Proposition~\ref{prop:compact:basic}(b).
\end{proof}

To conclude, we investigate the class of topological $*$-algebras whose
$*$-representation topology is generated by a single continuous seminorm. 
As the next lemma (which is modeled on \cite[10.1.7]{Palmer2}) reveals, 
these are precisely the algebras for which $\overline{PC^*}(A)$ is 
a $C^*$-algebra.

\begin{lemma} \label{lemma:GG}
Let $A \in \mathsf{T^*A}$. The following statements are equivalent:

\begin{list}{{\rm (\roman{enumi})}}
{\usecounter{enumi}\setlength{\labelwidth}{30pt}\setlength{\topsep}{-0pt}
\setlength{\itemsep}{-0pt}\setlength{\labelsep}{6pt}}

\item
For every $x\in A$,
\begin{equation}
\gamma(x)=\sup \{p(x)\mid p\in \mathcal{N}(A)\} = 
\{\|\pi(x)\| \mid \pi \in \mathcal{R}(A)\} < \infty,
\end{equation}
and $\gamma$ is continuous on $A$;

\item
$\mathcal{N}(A)$ has a largest element;

\item
$\overline{PC^*}(A)$ is  a $C^*$-algebra;

\item
$\mathcal{T}_A$ is defined by a single continuous linear space 
seminorm on $A$.
\end{list}
\end{lemma}

\pagebreak[2]

\begin{remark}
The seminorm in (iv) is not assumed {\em a priori}
to be submultiplicative or a $C^*$-seminorm.
\end{remark}

\begin{proof}
(i) $\Rightarrow$ (ii)
Since each $p\in\mathcal{N}(A)$ is a $C^*$-seminorms, so is $\gamma$, which
(being continuous) belongs to $\mathcal{N}(A)$. Thus, $\gamma$ is 
the largest element of $\mathcal{N}(A)$.

(ii) $\Rightarrow$ (iii)
Let $r\in \mathcal{N}(A)$ be the largest element. Then $A_r$ is a $C^*$-algebra, 
and for each $p\in \mathcal{N}(A)$ there is a morphism $A_r \rightarrow A_p$. 
Thus, $\overline{PC^*}(A)=\lim A_*=A_r$ is a $C^*$-algebra.

(iii) $\Rightarrow$ (iv)
Let $\| \cdot\|$ be the norm of the $C^*$-algebra $\overline{PC^*}(A)$. 
It is certainly a $C^*$-seminorm on $A$, and by definition,
it defines $\mathcal{T}_A$. In particular, $\|\cdot\|$
is continuous.

(iv) $\Rightarrow$ (i)
Let $\sigma$ be a seminorm that defines $\mathcal{T}_A$. Since multiplication
is jointly continuous in $PC^*(A)$ (which carries the quotient topology
$\mathcal{T}_A/A_R$), there is a constant $B>0$ such that 
$\sigma(xy)\leq\nolinebreak B \sigma(x)\sigma(y)$
for every $x,y\in A$. Thus, by replacing $\sigma$ with $B\sigma$ if necessary, 
we may assume that $B=1$ and $\sigma$ is submultiplicative. Let 
$\pi\colon A \rightarrow B(H) \in \mathcal{R}(A)$.  Pick $a\in A$.
Since $B(H)$ is a $C^*$-algebra,
\begin{align}
\|\pi(a)\|^2  = & \|\pi(a^* a)\| = r_{B(H)}(\pi(a^* a)) = 
\lim\limits_{n \rightarrow \infty} \|(\pi(a^* a))^n\|^{\frac 1 n}.\\
\intertext{By the continuity of $\pi$ with respect to $\mathcal{T}_A$, 
there is a constant $C_\pi$ such that
$\|\pi(x)\| \leq C_\pi \sigma(x)$ for every $x \in A$, so}
& \|(\pi(a^* a))^n\|^{\frac 1 n}  \leq 
\sqrt[n]{C_\pi} \sigma((a^* a)^n)^{\frac 1 n}
\leq \sqrt[n]{C_\pi} \sigma(a^* a).\\
\intertext{Therefore,}
\|\pi(a)\|^2  = & 
\lim\limits_{n \rightarrow \infty} \|(\pi(a^* a))^n\|^{\frac 1 n} \leq
\lim\limits_{n \rightarrow \infty} \sqrt[n]{C_\pi} \sigma(a^* a) = 
\sigma(a^* a) \leq D \sigma(a)^2,
\end{align}
where $D$ is a constant such that $\sigma(x^*) \leq A\sigma(x)$ for every
$x \in A$. Hence, $\gamma(a) \leq \sqrt D \sigma(a)$, which proves 
the second statement too, because $\sigma$ is assumed to be 
continuous on $A$. 
\end{proof}

Palmer \cite[10.1]{Palmer2} investigated discrete $*$-algebras that satisfied
the equivalent conditions of Lemma~\ref{lemma:GG}, and named them
{\em $G^*$-algebras}. Thus, we extend this terminology to 
topological $*$-algebras too, and denote by $\mathsf{G^*A}$ the full 
subcategory of $\mathsf{T^*A}$ formed by the (generalized)
$G^*$-algebras. For $A \in \mathsf{G^*A}$, we put $C^*(A)$ for 
$\overline{PC^*}(A)$, and call it the {\em enveloping $C^*$-algebra} of $A$.
The next Proposition is a consequence of Proposition~\ref{prop:TA-ref}
and Lemma~\ref{lemma:GG}; its restricted version, for discrete
$G^*$-algebras, appears in \cite[10.1.11]{Palmer2}.

\begin{proposition}
The category $\mathsf{C^*A}$ is a full reflective subcategory of 
$\mathsf{G^*A}$, and its reflector is given by $A \longmapsto C^*(A)$.
\qed
\end{proposition}

We once again return to the setting of section~\ref{sect:basic}.
Equip $\cat{X}=\mathsf{G^*A}$ with the 
$(\mathsf{Onto},\mathsf{Embed})$-factorization system, and provide
$\cat{A}=\mathsf{C^*A}$ with the 
$(\mathsf{Onto},\mathsf{Inj})$-factorization system
($\mathsf{Inj}$ stands for injections). Notice that
every injection of $C^*$-algebras is actually an embedding
(cf.~\cite[1.8.1]{DixmierC}), so the Kuratowski closure on
$C^*$-subalgebras of $C^*$-algebras is just the discrete closure $s$.
(This is not too surprising, though, because van Osdol \cite{Osdol} showed 
that the category of unital $C^*$-algebras and their unital $*$-homomorphisms
is monadic over  $\,\mathsf{Set}$.) Therefore, by 
Corollary~\ref{cor:compact:trivial}, we get:

\begin{theorem}
An algebra $A\in\mathsf{G^*A}$ is $c^*$-compact if and only if
$\alpha_A\colon A \rightarrow C^*(A)$ is surjective. \qed
\end{theorem}

%01234567890123456789012345678901234567890123456789012345678901234567890123456789

\section*{Acknowledgments}

I am deeply indebted to Professor Walter Tholen, my PhD supervisor, who
introduced me to the concept of closure operators, for his attention to 
this work and his helpful suggestions.

I am grateful to Professor Dikran Dikranjan, whose talk at CITA-2003 was
the main source of inspiration for this research, for his sharp and
constructive criticism of my early notes for the paper, and for
numerous valuable suggestions at the later stages.

I wish to thank Professor Salvador Hern\'andez, Professor Maria Manuel
Clementino, and Professor Bob Par\'e for the valuable discussions that 
were of great assistance in writing this paper.

\bibliography{notes,notes2,notes3}

\begin{samepage}

{\bigskip\bigskip\noindent 
FB 3 - Mathematik und Informatik\\
Universit\"at Bremen\\ 
Bibliothekstrasse 1 \\
28359 Bremen \\
Germany

\bigskip\noindent
{\em Current address:}\\
Department of Mathematics and Statistics\\
Dalhousie University\\
Halifax, B3H 3J5, Nova Scotia\\
Canada

\nopagebreak
\bigskip\noindent{\em e-mail: lukacs@mathstat.yorku.ca} }
\end{samepage}

\end{document}